\documentclass[12pt, a4paper]{article}

\usepackage{amssymb}
\usepackage{amsmath}

  \newtheorem{theorem}{Theorem}
 \newtheorem{definition}{Definition}
  \newtheorem{corollary}{Corollary}
    \newtheorem{lemma}{Lemma}    
    
        \newtheorem{question}{Question}

\begin{document}

 \title{Interpolation, extrapolation, Morrey spaces and local energy control for the Navier--Stokes equations}
\author{Pierre Gilles Lemari\'e--Rieusset\footnote{LaMME, Univ Evry, CNRS, Universit\'e Paris-Saclay, 91025, Evry, France; e-mail : pierregilles.lemarierieusset@univ-evry.fr}}
\date{}\maketitle

\begin{abstract}
Barker recently proved new weak-strong uniqueness results for the Navier--Stokes equations based on  a criterion involving Besov spaces and a proof through interpolation between Besov-H\"older spaces and $L^2$. We improve slightly his results by considering Besov-Morrey spaces and interpolation between Besov-Morrey spaces and $L^2_{\rm uloc}$.
\end{abstract}
 
\noindent{\bf Keywords : } Navier--Stokes equations, Morrey spaces, Besov spaces, uniformly locally square integrable functions,  weak-strong uniqueness\\

\noindent{\bf AMS classification : } 35K55, 35Q30, 76D05.

\section{The Navier--Stokes equations}
 Let $\vec u_0$ a divergence-free vector field on $\mathbb{R}^3$. We shall consider weak  solutions to the Cauchy initial value problem for the Navier--Stokes equations which satisfy energy estimates.

The differential  Navier--Stokes equations read as
 \[ \partial_t \vec u+\vec u.\vec\nabla \vec u=\Delta \vec u-\vec\nabla p\]
 \[ \text{ div }\vec u=0\]
 \[ \vec u(0,.)=\vec u_0\]

Under reasonable assumptions, the problem is equivalent to the following integro-differential problem : 
\[ \vec u=e^{t\Delta}\vec u_0-B(\vec u,\vec u)(t,x)\] where
\begin{equation} \label{bilin} B(\vec u,\vec v)=\int_0^t e^{(t-s)\Delta} \mathbb{P} \text{ div }(\vec u\otimes\vec v) \, ds\end{equation} and $\mathbb{P}$ is the Leray projection operator.  (See \cite{LEM1, LEM2}  for details).

\subsection*{Weak Leray solutions for the Navier--Stokes equations}

  When $\vec u_0\in L^2$,                                                                                                                                                                                                                                                                                                                                                                                                                                                                  
 Leray proved existence of solutions $\vec u$ on $(0,+\infty)\times\mathbb{R}^3$ such that :
 \begin{itemize}
 \item $\vec u\in L^\infty_t L^2_x\cap L^2_t \dot H^1_x$
 \item $\lim_{t\rightarrow 0^+} \|\vec u(t,.)-\vec u_0\|_2=0$
 \item we have the Leray energy inequality
\begin{equation}\label{eqler} \|\vec u(t,.)\|_2^2+2\int_0^t \|\vec\nabla\otimes\vec u\|_2^2\, ds\leq \|\vec u_0\|_2^2\end{equation}
 \end{itemize}
 Such solutions are weak solutions : the derivatives in the Navier--Stokes solutions are taken in the sense of distributions. Those solutions (that satisfy the energy inequality (\ref{eqler})) are called  \textit{Leray weak solutions}.\\

  When $\vec u_0\in L^2$,                                                                                                                                                                                                                                                                                                                                                                                                                                                                  
 Leray's proof of   existence of solutions  \cite{LER}  is based on mollification, energy estimates and compactness arguments :
 \begin{itemize}
 \item he solves $$\partial_t \vec u_\epsilon+(\varphi_\epsilon*\vec u_\epsilon).\vec\nabla \vec u_\epsilon=\Delta\vec u_\epsilon -\vec\nabla p_\epsilon$$ with $\text{ div }\vec u_\epsilon=0$ and $\vec u_\epsilon(0,.)=\vec u_0$. Here, $\varphi\in \mathcal{D}$, $\int\varphi\, dx=1$ and $\varphi_\epsilon(x)=\frac 1{\epsilon^3} \varphi(\frac x \epsilon)$.
  \item  the solution holds on an interval $(0,T_\epsilon)$ where $T_\epsilon$ depends on $\epsilon$ and on $\|\vec u_0\|_2$  and we have 
the equality $$ \|\vec u_\epsilon(t,.)\|_2^2+2\int_0^t \|\vec\nabla\otimes\vec u_\epsilon\|_2^2\, ds= \|\vec u_0\|_2^2$$
\item the solution is then global; moreover by Rellich theorem, we find a subsequence that converges strongly in $(L^2_t L^2_x)_{loc}$ to a Leray solution $\vec u$
 \end{itemize}
   Such solutions (i. e. obtained by this mollification/extraction process) will be called in the following \textit{restricted Leray weak solutions}.\\

  Restricted Leray solutions satisfy the Leray energy inequality which takes into account the energy on the whole space.
  But they enjoy as well a pointwise inequality property : for a non-negative locally finite measure $\mu$ we have
\begin{equation}\label{eneq} \partial_t(\vert\vec u\vert^2)+ 2 \vert\vec\nabla\otimes\vec u\vert^2=\Delta(\vert\vec u\vert^2)-\text{ div }((2p+\vert\vec u\vert^2)\vec u)-\mu\end{equation}
Leray solutions that enjoy the pointwise energy inequality are called  \textit{suitable Leray solutions} \cite{CKN}.

\subsection*{Local weak Leray solutions}
 The pointwise energy inequality allows one  \cite{LEM0, LEM1} to develop a theory of weak solutions with infinite energy. Consider $\vec u_0$ a divergence-free vector field that is uniformly locally square integrable :
 \[ \sup_{x_0\in\mathbb{R}^3} \int_{\vert x-x_0\vert<1} \vert\vec u_0(x)\vert^2\, dx<+\infty.\]
  
A  \textit{local Leray solution} on $(0,T)\times\mathbb{R}^3$ is a solution such that
   \begin{itemize}
 \item $\vec u\in L^\infty_t (L^2_{\rm uloc})\cap (L^2_t \dot H^1_x)_{\rm uloc}$
 \item for all compact subset $K$ of $\mathbb{R}^3$,  $\lim_{t\rightarrow 0^+} \int_K\vert \vec u(t,.)-\vec u_0\vert^2\, dx=0$
 \item we have the pointwise  energy inequality   (\ref{eneq}).
 \end{itemize}
  Local in time existence of   {restricted local Leray solutions} has been proved for a positive $T$ that depends only on $\|\vec u_0\|_{L^2_{\rm uloc}}$ (see section \ref{local}).

 \section{The Prodi--Serrin criterion for weak-strong uniqueness} Based on a compactness criterion, the proof of existence of Leray solutions does not provide any clue on the would-be uniqueness of the solution to the Cauchy initial value problem.

 A classical case of uniqueness of Leray weak solutions is
  {Serrin's criterion for weak-strong uniqueness} \cite{PRO, SER}.
If $\vec u_0\in L^2$ and if the Navier-Stokes equations has a solution $\vec u$ on $(0,T)$  such that $$\vec u\in X_T= L^p_t L^q_x\text{ with } \frac 2 p+\frac 3 q=1 \text{ and } 2\leq p<+\infty$$ then if $\vec v$ is a Leray solution we have $\vec u=\vec v$ on $(0,T)$.

The proof of the criterion is based on the fact that if $\vec v$ is a Leray solution and if $\vec u$ is the mild solution with $\|\vec u\|_{X_T}<+\infty$, then
the difference $\vec w=\vec u-\vec v$ satisfies a Gronwall estimate :
$$ \|\vec w(t,.)\|_2^2 +2 \int_0^t \|\vec \nabla\otimes\vec w\|_2^2\, ds\leq 2 \int_0^T\vert \int \vec u.(\vec w.\vec\nabla\vec w)\, dx\vert\, ds.$$
We have (for $\frac 2 p+\frac 3 q=1$)
$$ \|\vec u\otimes \vec w\|_2\leq C \|\vec u\|_q \|\vec w\|_2^{\frac 2 p} \|\vec\nabla\otimes \vec w\|_2^{\frac 3 q}$$ so that
\[ \vert \int \vec u.(\vec w.\vec\nabla \vec w)\ dx\vert\leq \|\vec u\|_q \|\vec w\|_2^{\frac 2 p} \|\vec\nabla\otimes \vec w\|_2^{1+\frac 3 q} = \|\vec u\|_q (\|\vec w\|^2)^{\frac 1 p} (\|\vec\nabla\otimes w\|^2)^{1-\frac 1 p}. \] Thus, we find
$$2 \int_0^T\vert \int \vec u.(\vec w.\vec\nabla\vec w)\, dx\vert\, ds\leq C \int_0^t \|\vec u\|_q^p \|\vec w\|_2^2\, ds+ \int_0^t \|\vec\nabla\otimes\vec w\|_2^2\, ds.$$
The facts that $\|\vec w(0,.)\|_2^2=0$ and $ \|\vec w(t,.)\|_2^2 \leq C \int_0^t \|\vec u\|_q^p\, \|\vec w\|_2^2\, ds$ then gives $\vec w=0$ on $(0,T)$.

We may comment a little further in the case $2<p<+\infty$. In that case, the bilinear operator $B$ (given by (\ref{bilin})) is bounded on $X_T=L^p_t L^q_x$. Thus, we find that 
 the existence of $T>0$ and of a solution in $L^p L^q$ with $\frac 2 p+\frac 3 q$ is equivalent to the existence of $T'$ such that $e^{t\Delta}\vec u_0\in L^p L^q$ on $(0,T')$ (and on $(0,+\infty)$, since $\vec u_0\in L^2$). Using the thermic characterization of Besov spaces, we can see that this is equivalent  with $$\vec u_0\in \dot B^{-1+\frac 3 q}_{q,p}.$$ Thus, the initial value is not only in $L^2$, but it must belong as well to a Besov space with a better regularity than provided by the embedding $$L^2\subset \dot B^{-\frac 3 2+\frac 3 q}_{q,2} \subset \dot B^{-\frac 3 2+\frac 3 q}_{q,p}.$$

\section{ The Koch and Tataru theorem and T. Barker's question} We may now wonder how to generalize the Prodi--Serrin criterion. It means : given $\vec u_0\in L^2$ and weak Leray solutions associated to $\vec u_0$, find a space $\mathbb{X}$ (as large as possible) such that if moreover $\vec u_0\in \mathbb{X}$ then we have a solution $\vec u\in X_T$ for some space $X_T$ of functions on $(0,T)\times \mathbb{R}^3$  and such that the existence of a solution in $X_T$ implies that any other weak Leray solution is equal to this solution $\vec u$ for $0<t<T$.

\subsection*{The space ${\rm BMO}^{-1}$}  First of all, we precise which kind of space $\mathbb{X}$ we are going to study.  The idea is to look at an initial value which generates a  solution in  some uniqueness class (where uniqueness holds for small  solutions). The setting where to construct such solutions is the setting of mild solutions, as introduced by Kato \cite{KAT}  : mild solutions are constructed by the Banach contraction principle.

Due to the symmetries of the equations (if $\vec u$ is a solution for initial value $\vec u_0$, then $\lambda \vec u(\lambda^2 t, \lambda(x-x_0))$ is a solution for the initial value $\lambda\vec u_0(\lambda(x-x_0))$), we look for spaces with norms invariant through the transforms $\vec u_0(.)\mapsto \lambda \vec u(\lambda(.-x_0))$ (for $\lambda>0$). Moreover, in order to be able to define
\[ B(\vec u,\vec v)=\int_0^t e^{(t-s)\Delta} \mathbb{P} \text{ div }(\vec u\otimes\vec v) \, ds\]  at least for $\vec u=\vec v=e^{t\Delta}\vec u_0$ (first step of the Picard iteration to find a fixed-point to $\vec u=e^{t\Delta}\vec u_0-B(\vec u,\vec u)$), we ask that $\int_{[0,1]\times B(0,1)} \vert e^{s\Delta}\vec u_0(y)\vert^2\, ds\, dy<+\infty$.\\

Thus, we are lead to introduce the space $\mathbb{X}$ of distributions  $v$ such that
\begin{itemize}
\item  $\sup_{t>0} \sqrt t \|e^{t\Delta} v \|_\infty <+\infty$
\item  $\sup_{0<t, x_0\in \mathbb{R}^3}  t^{-3/2} \int_0^t \int_{B(x_0,\sqrt t)} \vert  e^{s\Delta} v(y)\vert^2\, dy\, ds)^{1/2}$
\end{itemize}

This space $\mathbb{X}$ has been identified by Koch and Tataru \cite{KOT} :  this is the Triebel--Lizorkin space $\dot F^{-1}_{\infty,2}$, or equivalently the space ${\rm BMO}^{-1}=\sqrt{-\Delta} \, {\rm BMO}$. Moreover, they proved the following theorem :

\begin{theorem} \label{theoKT}
 For $0<T\leq\infty$, define
\[ \|\vec u\|_{X_T}=\]\[\sup_{0<t<T} \sqrt t \|\vec u(t,.\|_\infty + \sup_{0<t<T, x_0\in \mathbb{R}^3}  (t^{-3/2} \int_0^t \int_{B(x_0,\sqrt t)} \vert \vec u(s,y)\vert^2\, dy\, ds)^{1/2}
\]

There exists $C_0$ (which does not depend on $T$)  such that if $T\in (0,+\infty]$, if $\vec u$ and $\vec v$ are defined on $(0,T)\times\mathbb{R}^3$ and if
\[ B(\vec u,\vec v)=\int_0^t e^{(t-s)\Delta} \mathbb{P} \text{ div }(\vec u\otimes\vec v) \, ds\] 
then
\[ \|B(\vec u,\vec v)\|_{X_T}\leq C_0 \|\vec u\|_{X_T}\|\vec v\|_{X_T}.\]
\end{theorem}

\begin{corollary}  If $\|e^{t\Delta}\vec u_0\|_{X_T}<\frac 1{4C_0}$, then the integral Navier--Stokes equations have a solution on $(0,T)$ such that $\|\vec u\|_{X_T}\leq 2\| e^{t\Delta}\vec u_0\|_{X_T}$.
 
    This is the  {unique} solution such that $\|\vec u\|_{X_T}\leq \frac 1 {2C_0}$.
\end{corollary}
 
 A special case of initial data that leads to a solution in some $X_T$ is given by the subspace ${\rm VMO}^{-1}$ of ${\rm BMO}^{-1}$.
 
 \begin{definition} 
 ${\rm VMO}^{-1}$ is the closure of compactly supported    functions in ${\rm BMO}^{-1}$.
 \end{definition}   If $\vec u_0\in {\rm VMO}^{-1}$, then $\lim_{T\rightarrow 0} \|e^{t\Delta}\vec u_0\|_{X_T}=0$. 
Remark that we have the embedding $\dot B^{-1+\frac 3 q}_{q,p} \subset {\rm VMO}^{-1}$ for $2<p<+\infty$ and $\frac 2 p+\frac 3 q=1$. As a matter of fact, we may consider ${\rm VMO}^{-1}$ as a limit case for the scale of spaces $\dot B^{-1+\frac 3 q}_{q,p}$. Thus, Barker \cite{BAR} raised the following question :

\begin{question} If $\vec u_0$ belongs to $L^2\cap {\rm VMO}^{-1}$, does there exists a positive time $T$ such that every weak Leray solution of the Cauchy problem for the Navier--Stokes equations with $\vec u_0$ as initial value coincide with the mild solution in $X_T$ ?
\end{question}

If $\vec u_0\in L^2\cap {\rm VMO}^{-1}$ and if $\|e^{t\Delta}\vec u_0\|_{X_T}\leq \frac 1 {4C_0}$, then if $\vec u$ is a restricted Leray solution of the Navier--Stokes solutions with initial value $\vec u_0$, then $\|\vec v\|_{X_T}\leq   2 \|e^{t\Delta}\vec u_0\|_{X_T}$. In particular, we have uniqueness of \emph{restricted Leray solutions } on $(0,T)$. 

As a matter of fact, this proof of local uniqueness of restricted  weak Leray solutions holds for a slightly more general class :

\begin{definition} 
 ${\rm BMO}^{-1}_0$ is the space of distributions $u_0$  in ${\rm BMO}^{-1}$ such that  $$\lim_{T\rightarrow 0} \|e^{t\Delta}  u_0\|_{X_T}=0.$$
 \end{definition}

  \textit{In the following, we will focus on the hypothesis $\vec u_0\in L^2\cap {\rm BMO}^{-1}_0$ and on the issue of uniqueness for Leray solutions.}
 
\section{The limiting case}
  Up to now, we don't know how to prove local uniqueness of the Leray solutions when the initial value $\vec u_0$ belongs to $L^2\cap {\rm BMO}^{-1}_0$.  What we know for sure is that the mild solution $\vec u$ in $X_T$ belongs to $L^\infty((\epsilon,T)\times \mathbb{R}^3)$ for every positive $\epsilon  \in (0,T)$. Moreover, $\vec u$ is a weak  Leray solution and  for every other   Leray solution $\vec v$ and for $\epsilon>0$, we have 
    $$\partial_t(\vec u.\vec v)=\vec u.\partial_t\vec v+\partial_t\vec u.\vec v$$
    which gives for $0<\epsilon<t<T$
    \begin{equation*}\begin{split} \int \vec u(t,x).\vec v(t,x)\, dx=\int \vec u(\epsilon,x).\vec v(\epsilon,x) \, dx& -2\int_\epsilon^t\int  \vec\nabla\otimes\vec u.\vec \nabla\otimes\vec v\, dx\, ds\\& -\int_\epsilon^t \int \vec u.(\vec v.\vec\nabla \vec v)+\vec v.(\vec u.\vec\nabla\vec u)\, dx\, ds
  \end{split}\end{equation*} From 
  $$ \int \vec u.(\vec v.\vec\nabla\vec v)+\vec v.(\vec u.\vec\nabla\vec u)\, dx \!=\!\int \vec u.(\vec v.\vec\nabla ( \vec  v-\vec u))+(\vec v-\vec u).(\vec u.\vec\nabla\vec u)\, dx\! =\! \int\vec u\big( (\vec v-\vec u).\vec\nabla( \vec  v-\vec u)\big) \, dx$$ and letting $\epsilon$ go to $0$, we get
      \begin{equation*}\begin{split} \int \vec u(t,x).\vec v(t,x)\, dx = \|\vec u_0\|_2^2& -2\int_0^t\int  \vec\nabla\otimes\vec u.\vec \nabla\otimes\vec v\, dx\, ds\\&  -\lim_{\epsilon\rightarrow 0}\int_\epsilon^t \int\vec u\big( (\vec v-\vec u).\vec\nabla( \vec  v-\vec u)\big) \, dx  \, ds
  \end{split}\end{equation*} and (letting $\vec v=\vec u$)
  $$ \|\vec u(t,.)\|_2^2=\|\vec u_0\|_2^2-2 \int_0^t \|\vec\nabla\otimes \vec u\|_2^2 \, ds.$$
  Combining those two equalities with the Leray energy inequality for $\vec v$
  \begin{equation*}  \|\vec v(t,.)\|_2^2+2\int_0^t \|\vec\nabla\otimes\vec v\|_2^2\, ds\leq \|\vec u_0\|_2^2\end{equation*}
  we get the following inequality for $\vec w=\vec u-\vec v$ :
  \begin{equation}  \|\vec w(t,.)\|_2^2 +2 \int_0^t \|\vec \nabla\otimes\vec w\|_2^2\, ds\leq 2 \int_0^T\vert \int \vec u.(\vec w.\vec\nabla\vec w)\, dx\vert\, ds.
  \end{equation}

  As a matter of fact,  the key ingredient in Prodi--Serrin's criterion is the estimation of the integral
  $$ I(\vec u,\vec w)= \int_0^t \left\vert \int \vec u(.\vec w.\vec\nabla\vec w) \, dx\right\vert\, ds$$ but, if $\vec w=\vec v-\vec u$ with $\vec v$ a Leray solution and $\vec u$ the mild solution in $X_T$, we don't even know whether $I(\vec u,\vec w)$ is finite.
  
  In the limiting case of Prodi and Serrin,  (for $p=2$ and $q=+\infty$),  we write
  $$ \|\vec w(t,.)\|_2^2 + 2\int_0^t \|\vec\nabla\otimes\vec w\|_2^2\, ds\leq 2 \int_0^t \int_0^t \left\vert \int \vec u(.\vec w.\vec\nabla\vec w) \, dx\right\vert\, ds \leq 2 \int_0^t \|\vec u\otimes\vec w\|_2 \|\vec\nabla\otimes\vec w\|_2\, ds$$
  and get 
\begin{equation} \label{gronw} \|\vec w(t,.)\|_2^2  \leq  \int_0^t \|\vec u\|_\infty^2 \|\vec w\|_2^2\, ds.\end{equation}

  Of course, we may conclude under 
  the assumption that $\vec u\in L^2_t L^\infty_x$. Actually, we shall not be interested in measurabilty issues for functions with values in a non-separable space such as $L^\infty$ (i.e. in Bochner measurability for instance), as we are dealing with locally integrable functions for the Lebesgue measure $dt\, dx$ on $(0,T)\times\mathbb{R}^3$. Thus, for almost every $t$ the quantity $\|\vec u(t,.)\|_\infty$ will be wel-defined as a measurable function of $t$, and $\vec u\in L^2_t L^\infty_x$ will simply mean that $\int_0^T \|\vec u(t,.)\|_{L^\infty(dx)}^2\, dt<+\infty$. (See \cite{LEM1} for details.)
  
  If $\vec u\in L^2_t L^\infty_x$ on $(0,T)\times\mathbb{R}^3$, the Prodi--Serrin criterion proves that every weak Leray solution $\vec v$ on $(0,T)$ is equal to the mild solution $\vec u$. (This has even be extended to the case $\vec u\in L^2_t {\rm BMO}_x$ by Kozono and Tanyuchi \cite{KOTA}).    But it is not  easy to translate the condition that $\vec u\in L^2_t L^\infty_x$ into an equivalent assumption on $\vec u_0$. The problem comes from the fact that the bilinear operator $B$ is not bounded on $L^2_t L^\infty_x$. 
  
  On the other hand, if we only assume $\vec u_0\in L^2\cap {\rm BMO}^{-1}_0$, we only know the inequality  $\|\vec u(t,.)\|_\infty\leq \|\vec u\|_{X_T} \frac {o(1)}{\sqrt t}$. We find an integrability issue near $t=0$. To check that this is actually an issue, consider the following example : take $\vec\omega$ a divergence-free vector field in the Schwartz class such that the Fourier transform of $\vec\omega$ is compactly supported in the annulus $1<\vert\xi\vert<2$; define
  $$\vec u_0(x)=\sum_{j=0}^{+\infty}  2^j\frac 1{\sqrt{1+j}}\vec \omega(2^jx)  ;
$$  
 we have $\vec u_0\in L^2$;    $\vec u_0 \in \dot B^{-1+\frac 3 q}_{q,p}\subset \dot B^{-1}_{\infty,p}$ (with $2<p<+\infty$ and $\frac 2 p+\frac 3 q=1$) but $\vec u_0\notin \dot B^{-1}_{\infty,2}$;  as the bilinear operator $B$ is bounded on $L^p_tL^q_x\cap L^2_t L^\infty_x$, if we assume that the mild solution $\vec u\in L^p_tL^q_x$ belongs to $L^2_t L^\infty_x$, we would  find that $e^{t\Delta}\vec u_0\in L^2_t L^\infty_x$, and thus $\vec u_0\in \dot B^{-1}_{\infty,2}$; thus, we have $\int_0^T \|\vec u\|_\infty^2\, dt=+\infty$.

 \section{Barker's theorem}   In this section, we  shall sketch the proof of Barker \cite{BAR}, as we shall extend it in Section \ref{mainsec} to the case of Besov--Morrey spaces.The main idea in the recent paper of Barker  is the following one : if we want to use only the inequality
$$  \|\vec u(t,.)\|_\infty\leq \|\vec u\|_{X_T} \frac {o(1)}{\sqrt t}$$  to deal with the Gronwall inequality (\ref{gronw}), we need to assume more than $\vec w\in L^\infty_t L^2_x$. Indeed, we have the easy following lemma :

\begin{lemma} \label{lemgron} Let $\delta>0$. Let $A$ and $B$  be   locally  bounded non-negative measurable functions on $(0,T]$ such that $$\lim_{t\rightarrow 0} t A(t)=0 \text{ and } \sup_{ 0<t<T } t^{-\delta} B(t)=0.$$ If we have moreover, for all $t\in (0,T]$,
$$ B(t)\leq \int_0^t A(s) B(s)\, ds$$ then $B=0$.
\end{lemma}
\noindent{\bf Proof : }We have 
$$ B(t)\leq \frac {t^\delta} \delta \sup_{0<s<t} s A(s) \sup_{0<\sigma<t} \sigma^{-\delta} B(\sigma)$$ so that $B=0$ on $(0,T_0]$ as long as $\sup_{0<s<T_0} sA(s) <\delta$. For $t>T_0$, we then write $B(t)\leq \sup_{T_0<s<T} A(s) \int_{T_0}^t B(s)\, ds$ and we find $B=0$.
$\diamond$

As $\vec w=\vec v-\vec u=(\vec v-e^{t\Delta}\vec u_0)-(\vec u-e^{t\Delta}\vec u_0)$, the extra information on $\|\vec w\|_2$  will be provided by the following lemma :
\begin{lemma} \label{lembar}  Let $\vec u_0$ be a divergence-free vector field with $\vec u_0\in L^2$and let $\vec v$ be a weak Leray solution of the Navier--Stokes equations with initial value $\vec u_0$.
If   moreover $$\vec u_0\in   [L^2, B^{-\gamma}_{\infty,\infty}]_{\theta,\infty} \text{ for some } -1<-\gamma<0\text{ and } 0<\theta<1$$
 then there exists  $\delta>0$  such that
 $$  \sup_{t>0} t^{-\delta}\|\vec v(t,.)-e^{t\Delta}\vec u_0\|_2<+\infty.$$
\end{lemma}

\noindent{\bf Proof : } Assume $\vec u_0\in   [L^2, B^{-\gamma}_{\infty,\infty}]_{\theta,\infty} $.  For $0<\epsilon<1$, split $\vec u_0$ in $\vec \alpha_\epsilon+\vec\beta_\epsilon$ with
 $$ \|\vec\alpha_\epsilon\|_2\leq C_0 \epsilon^{\theta}\text{ and } \|\vec \beta_\epsilon\|_{\dot B^{-\gamma}_{\infty,\infty}}\leq C_0 \epsilon^{\theta-1}$$ where $C_0$ does not depend on $\epsilon$ (but depends on $\vec u_0$).

 We have a solution $\vec U_\epsilon$ for the Navier-Stokes equations with initial value $\vec\beta_\epsilon$ 
 such that $\| \vec U_\epsilon(t,.)\|_\infty\leq C_1 t^{-\frac {\gamma}2}  \epsilon^{\theta-1}$ on an interval $(0,T_\epsilon)$  with $T_\epsilon^{\frac {1-\gamma}2} \epsilon^{\theta-1}= C_2$.
  Moreover $\vec\beta_\epsilon=\vec u_0-\vec\alpha_\epsilon\in L^2$ and we find that $\sup_{0<t<T_\epsilon} \|\vec U_\epsilon\|_2\leq C_3$ and that $\|\vec U_\epsilon-e^{t\Delta}\vec\beta_\epsilon\|_2\leq C_4 t^{(1-\gamma)/2} \epsilon^{\theta-1}$.
 
 Since $\vec U_\epsilon$ is in $L^2_t L^\infty_x$ on every bounded interval, 
  we get
      \begin{equation*}\begin{split} \int \vec U_\epsilon(t,x).\vec v(t,x)\, dx = \int \vec\beta_\epsilon(x).\vec u_0(x)\, dx& -2\int_0^t\int  \vec\nabla\otimes\vec U_\epsilon.\vec \nabla\otimes\vec v\, dx\, ds\\&  - \int_0^t \int\vec U_\epsilon\big( (\vec v-\vec U_\epsilon).\vec\nabla( \vec  v-\vec U_\epsilon)\big) \, dx  \, ds
  \end{split}\end{equation*} and    $$ \|\vec U_\epsilon(t,.)\|_2^2=\|\vec \beta_\epsilon\|_2^2-2 \int_0^t \|\vec\nabla\otimes \vec U_\epsilon\|_2^2 \, ds.$$
  
  Combining those two equalities with the Leray energy inequality for $\vec v$
  \begin{equation*}  \|\vec v(t,.)\|_2^2+2\int_0^t \|\vec\nabla\otimes\vec v\|_2^2\, ds\leq \|\vec u_0\|_2^2\end{equation*}
  we get the following inequality for  $\vec W_\epsilon=\vec v-\vec U_\epsilon$:
  \begin{equation*}  \|\vec W_\epsilon(t,.)\|_2^2 +2 \int_0^t \|\vec \nabla\otimes\vec W_\epsilon\|_2^2\, ds\leq \|\vec \alpha_\epsilon\|_2^2+ 2 \int_0^T\vert \int \vec U_\epsilon.(\vec W_\epsilon.\vec\nabla\vec W_\epsilon)\, dx\vert\, ds.
  \end{equation*}
Thus, we get
$$   \|\vec W_\epsilon(t,.)\|_2^2 \leq C_0^2 \epsilon^{2\theta} + C_1^2 \epsilon^{2(\theta-1)} \int_0^t s^{-\gamma} \|\vec W_\epsilon(s,.)\|_2^2 \, ds$$ so that
$$   \|\vec W_\epsilon(t,.)\|_2^2 \leq C_0^2 \epsilon^{2\theta}  e^{ C_1^2 \epsilon^{2(\theta-1)}\frac{t^{1-\gamma}}{1-\gamma}}.
$$

 Now, for $\tau<1$, take $\epsilon=\tau^\mu$ with $\frac{1-\gamma} 2+\mu(\theta-1)>0$. We find that, for $0<t<T_\epsilon$ with 
 $$ T_\epsilon=  C_2^{\frac 2{1-\gamma}}\epsilon^{\frac{2 (1-\theta)}{1-\gamma}}= C_2^{\frac 2{1-\gamma}} \tau^{\mu \frac{2 (1-\theta)}{1-\gamma}}\  [\text{ where }  \mu \frac{2 (1-\theta)}{1-\gamma}<1],$$ we have the inequality
\begin{equation*}\begin{split} \|\vec v-e^{t\Delta}\vec u_0\|_2\leq & \|\alpha_\epsilon\|_2+\|\vec W_\epsilon\|_2+ \|\vec U_\epsilon-e^{t\Delta} \vec \beta_\epsilon\|_2
\\ \leq & C_0\epsilon^\theta (1+e^{ \frac {C_1^2}{2(1-\gamma)}\frac {t^{1-\gamma}}{T_\epsilon^{1-\gamma}}}) +C_4 t^{\frac{1-\gamma}2}\epsilon^{\theta-1}
\end{split}\end{equation*}
 If $\tau$ is small enough, we have $\tau<T_\epsilon$ and we find
 $$ \|\vec v(\tau,.)-e^{\tau\Delta}\vec u_0\|_2\leq    C_0\tau^{\mu\theta} (1+e^{ \frac {C_1^2}{2(1-\gamma }} )+C_4 \tau^{\frac{1-\gamma}2+\mu(\theta-1)} $$
 The lemma is proved.
$\diamond$
Barker's theorem then reads as :
   \begin{theorem} \label{theobar}  Let $\vec u_0$ be a divergence-free vector field with $\vec u_0\in L^2$and let $\vec v$ be a weak Leray solution of the Navier--Stokes equations with initial value $\vec u_0$.
If   moreover $$\vec u_0\in  {\rm BMO}^{-1}_0\cap  \dot B^{-s}_{q,\infty} \text{ with } 3<q<+\infty \text{ and }-s>-1+\frac 2 q$$
 then  there exists $T>0$ such that if   $\vec v$ is a  Leray solution and if $\vec u$ is the mild solution with $\|\vec u\|_{X_T}<+\infty$, then $\vec u=\vec v$ on $(0,T)$.
 \end{theorem}
 
 \noindent{\bf Remark} We have the embeddings $L^2\subset \dot B^{-\frac 3 2+\frac 3 q}_{q,2}  \subset \dot B^{-\frac 3 2+\frac 3 q}_{q,\infty}$, so that the information conveyed by the hypothesis $\vec u_0\in \dot B^{-s}_{q,\infty}$ is  interesting only for the high frequencies of $\vec u_0$. Moreover the embeddings $L^2 = \dot B^0_{2,2}$ and ${\rm BMO}^{-1}\subset \dot B^{-1}_{\infty,\infty}$ gives that $$ L^2\cap   {\rm BMO}^{-1}_0\subset \dot B^{-1+\frac 2 q}_{q,q}\subset \dot B^{-1+\frac 2 q}_{q,\infty}.$$ Thus, the information conveyed by the hypothesis $\vec u_0\in \dot B^{-s}_{q,\infty}$ is   not contained in the assumption $\vec u_0\in L^2\cap   {\rm BMO}^{-1}_0$. Finally, if $-s>-1+\frac 3 q=\frac 2 p$, then we have  $$ L^2\cap   {\rm BMO}^{-1}_0\cap \dot B^{-s}_{q,\infty}\subset \dot B^{-1+\frac 2 q}_{q,q} \cap \dot B^{-s}_{q,\infty}\subset \dot B^{-1+\frac 3 q}_{q,p}.$$ Thus, the theorem is interesting only in the range $-1+\frac 2 q<-s \leq -1+\frac 3 q$, which corresponds to the gap between  $L^2\cap   {\rm BMO}^{-1}_0$. (where loca`luniqueness is conjectured to hold) and  $\dot B^{-1+\frac 3 q}_{q,p}$ (for which the Prodi--Serrin criterion shows that local uniqueness holds). 
 
 A  further remark is that we have the embedding $\dot B^{-1+\frac 3 q}_{q,\infty}\subset {\rm BMO}^{-1}$, so that we have a Prodi-Serrin criterion with $\vec u\in L^p_t L^q_x$ (with  $\frac 2 p+\frac 3 q=1$) replaced with  \[\sup_{0<t<T} t^{\frac 1 p} \|\vec u\|_q<+\infty \text{ and }\lim_{t\rightarrow 0} t^{\frac 1 p} \|\vec u(t,.)\|_q=0 .\]

\noindent{\bf Proof :}

The first step is the use of interpolation inequalities in order to be able to check that  $\vec u_0$  fulfills the assuptions of Lemma \ref{lembar}.
\begin{itemize}
\item   Since  $\vec u_0\in L^2=\dot B^0_{2,2}$ and $\vec u_0\in {\rm BMO}^{-1}\subset \dot B^{-1}_{\infty,\infty}$, we have $\vec u_0\in \dot B^{-1+\frac 2 q}_{q,q}$.
\item Since $\vec u_0\in \dot B^{-1+\frac 2 q}_{q,q}\cap \dot B^{-s}_{q,\infty}$, we have for $-1+\frac 2 q< -\sigma<-s$, $\vec u_0\in \dot B^{-\sigma}_{q,1}$.
\item  The Besov space $ \dot B^{-\sigma}_{q,1}$  is embedded in the Sobolev space $ \dot W^{-\sigma,q}$.
\item For $q<r<\infty$, we have $\dot W^{-\sigma,q}=[L^2,\dot W^{-\delta,r}]_{[\theta]}$ with $\frac 1 q=\frac {1-\theta}2+\frac\theta r$ and $-s=-\theta\delta$ (complex inteerpolation)
\item  Since $\dot W^{-\delta,r}\subset \dot B^{-1+\frac 3 r}_{\infty,\infty}$, we have  $[L^2,\dot W^{-\delta,r}]_{[\theta]}\subset [L^2,\dot B^{-\delta-\frac 3 r}_{\infty,\infty}]_{\theta,\infty}$ with 
$$ -\delta= -\frac s \theta= -s \frac{\frac 1 2-\frac 1 r}{\frac 1 2-\frac 1 q} =- \frac{(-s)}{-1+\frac 2 q}+O(\frac 1 r)$$ so that $-\delta-\frac 3 r>-1$ for $r$ large enough
\end{itemize}

We may now end the proof :
recall that if $\vec v$ is a Leray solution and if $\vec u$ is the mild solution with $\|\vec u\|_{X_T}<+\infty$, then
the difference $\vec w=\vec u-\vec v$ satisfies a Gronwall estimate :
\begin{equation*} \label{gronw} \|\vec w(t,.)\|_2^2  \leq  \int_0^t \|\vec u\|_\infty^2 \|\vec w\|_2^2\, ds.\end{equation*}

By Lemma \ref{lembar}, we have $ \|\vec v(t,.)-e^{t\Delta}\vec u_0\|_2=O(t^\delta)$ and $\|\vec  u(t,.)-e^{t\Delta}\vec u_0\|_2=O(t^\delta)$ for some positive $\delta$. On the other hand, we know that $\|\vec u(t,.)\|_\infty= o(\frac 1 {\sqrt t})$.  Using Lemma \ref{lemgron}, we find that $\vec w=0$, and $\vec v=\vec u$.$\diamond$

 \section{The Prodi--Serrin criterion for Besov--Morrey spaces}   Morrey spaces provide a natural tool for extending the Prodi--Serrin criterion. 
 
 \begin{definition}
 For $1< r\leq q<+\infty$, we define the Morrey space $\dot M^{r,q}$ as the space of Lebesgue measurable functions $f$  on $\mathbb{R}^3$ such that
 $$  \sup_{R>0, x_0\in\mathbb{R}^3} R^{\frac 3 q-\frac 3 r} (\int_{B(x_0,R)} \vert f(x)\vert^r\, dx)^{1/r}=\| f\|_{\dot M^{r,q}}<+\infty.$$
 Similarly, the space $\dot M^{1,q}$ is the space of locally finite Borelian (signed) measure $\mu$ such that 
  $$  \sup_{R>0, x_0\in\mathbb{R}^3} R^{\frac 3 q-  3 } (\int_{B(x_0,R)}  d\vert \mu\vert) =\|  \mu\|_{\dot M^{1,q}}<+\infty.$$
 Remark : For absolutely continuous measures $d\mu = f\, dx$ with $f\in L^1_{\rm loc}$, we have
   $$  \|\mu\|_{\dot M^{1,q}}= \sup_{R>0, x_0\in\mathbb{R}^3} R^{\frac 3 q-  3 } (\int_{B(x_0,R)}  \vert f(x)\vert \, dx) .$$
 \end{definition}
 
 The key inequality in the proof of the Prodi--Serrin criterion was the inequality (for all $w\in H^1$)
 $$ \|uw\|_2\leq C  \|u\|_q \|w\|_2^{\frac 2 p} \|\vec\nabla w\|_2^{\frac 3 q}$$
 with $\frac 2 p+\frac 3 q=1$ and $3<q<+\infty$. If we want to replace this inequality by a more general inequality
 $$ \|uw\|_2\leq  N(u) \|w\|_2^{\frac 2 p} \|\vec\nabla w\|_2^{\frac 3 q}$$ (again with $\frac 2 p+\frac 3 q=1$ and $3<q<+\infty$), then we proved in \cite{LEMM} that the existence of a finite $N(u)$ is equivalent to the fact that $u\in \dot M^{2,q}$, and moreover that
 $N(u) \approx \|u\|_{\dot M^{2,q}}$.

 This leads to the following  easy extension of the Prodi-Serrin criterion :

 \begin{theorem}
If $\vec u_0\in L^2$ and if the Navier-Stokes equations has a solution $\vec u$ such that $$\vec u\in L^p_t \dot M^{2,q}_x\text{ with } \frac 2 p+\frac 3 q=1\text{ and } 3<q<+\infty$$ then if $\vec v$ is a Leray solution we have $\vec u=\vec v$ on $(0,T)$.
\end{theorem}
 
 If $3< q<+\infty$, the existence of $T>0$ and of a solution in $L^p \dot M^{2,q}$ with $\frac 2 p+\frac 3 q$ is equivalent to the existence of $T'$ such that $e^{t\Delta}\vec u_0\in L^p \dot M^{2;q}$ on $(0,T')$ (and on $(0,+\infty)$, since $\vec u_0\in L^2$), thus with $$\vec u_0\in \dot B^{-1+\frac 3 q}_{\dot M^{2,q},p}.$$
 This Besov--Morrey space has been introduced  in 1994. by Kozono and Yamazaki \cite{KOY}. It is easy to check that,  for $2<p<+\infty$ and $\frac 2 p+\frac 3 q=1$, we have the inequality
 $$ \|e^{t\Delta} u_0\|_{X_T}\leq C_q (\int_0^T \| e^{t\Delta} u_0\|_{\dot M^{2,q}}^p\, dt)^{1:p}$$
 so that we have the embedding $\dot B^{-1+\frac 3 q}_{\dot M^{2,q},p} \subset {\rm BMO}^{-1}_0$ for $2<p<+\infty$ and $\frac 2 p+\frac 3 q=1$.

 \section{Barker's theorem and Besov-Morrey spaces} \label{mainsec}
 We shall extend Barker's theorem. 
\begin{theorem} \label{main}
If \begin{itemize}
\item $\vec u_0\in L^2\cap {\rm BMO}^{-1}_0$ 
\item  $3<q<+\infty$, $-s>-1+\frac 2 q$  and $\vec u_0\in \dot B^{-s}_{\dot M^{1,q},\infty}$
\end{itemize}
 then there exists $T>0$ such that if   $\vec v$ is a suitable Leray solution and if $\vec u$ is the mild solution with $\|\vec u\|_{X_T}<+\infty$, then $\vec u=\vec v$ on $(0,T)$.
 \end{theorem}
 
 \noindent{\bf Proof : }
 As we shall see, the proof is very similar to Barker's proof for Theorem \ref{theobar} \cite{BAR}. However, we shall meet some technical issues.

We sketch the proof :
\begin{itemize}
\item $\vec u_0\in L^2=\dot B^0_{2,2}$ and $\vec u_0\in {\rm BMO}^{-1}\subset \dot B^{-1}_{\infty,\infty}$, thus $\vec u_0\in \dot B^{-1+\frac 2 q}_{q,q}$
\item $\vec u_0\in \dot B^{-1+\frac 2 q}_{q,q}\cap \dot B^{-s}_{\dot M^{1,q},\infty}$ thus  $\vec u_0\in \dot B^{-\sigma}_{\dot M^{p,q},\infty}$ for $1<p<q$, $\frac 1 p =1-\theta+\frac \theta q$ and $-\sigma= -s(1-\theta)+\theta(-1+\frac 2 q)> -1+\frac 2 q$. We shall take $p>2$.
\item as $p<q$, we have   $\dot B^{-1+\frac 2 q}_{q,q}\ \subset \dot B^{1+\frac 2 q}_{q,p}$. Thus,   for $-1+\frac 2 q< -\gamma<-\sigma$, $$ \dot B^{-1+\frac 2 q}_{q,q}\cap  \dot B^{-\sigma}_{\dot M^{p,q},\infty}\ \subset \dot B^{-\gamma}_{\dot M^{p,q},1}\subset \dot W^{-\gamma,\dot M^{p,q}}.$$ 
\end{itemize}
  
We now encounter our first problem. We can no longer write $W^{-\gamma,\dot M^{p,q}}$ as a subspace of an interpolate space between $L^2$ and $\dot B^{-1+\delta}_{\infty,\infty}$.  More precisely, let us assume  $\dot B^{-\gamma}_{\dot M^{p,q},1}\subset [L^2, \dot B^{-1+\delta}_{\infty,\infty}]_{\theta,\infty}$; by homogeneity of the norms, we must have  $-\gamma-\frac 3 q=-(1-\theta)\frac 3 2 -\theta (1-\delta)$. We have $$ [L^2, \dot B^{-1+\delta}_{\infty,\infty}]_{\theta,\infty} \subset \dot B^{-\theta(1-\delta)}_{L^{r,\infty},\infty}$$ with $r=\frac 2 {1-\theta}$. In particular, for $u\in  [L^2, \dot B^{-1+\delta}_{\infty,\infty}]_{\theta,\infty} $, we have that $e^\Delta e^{t\partial_3^2} u$ goes to $0$ in $\mathcal{S}'$ when $t$ goes to $+\infty$. But if $3p\leq 2q$, if $u$ depends only on $(x_1,x_2)$ and not on $x_3$, and if $u\in \dot B^{-\gamma}_{\dot M^{p,\frac{2q}3},1}(\mathbb{R}^2)$, then $u\in \dot B^{-\gamma}_{\dot M^{p,q},1}(\mathbb{R}^3)$ and $e^{\Delta} e^{t\partial_3^2} u=e^\Delta u$. Thus, we have a contradiction.

We better  use complex interpolation and write that
$$ W^{-\gamma,\dot M^{p,q}} =[ \dot M^{2,\frac 2 p q}, \dot W^{-\rho,\dot M^{r,\frac r p q}}]^{[\theta]}$$ for $r>p$, $\frac {1-\theta}2+\frac\theta r=\frac 1 p$, $\gamma=\theta \rho$. (For interpolation of Morrey spaces, see \cite{LEMMOR1, LEMMOR2})

Then, we remark that $\dot M^{2,\frac 2 p q}\subset L^2_{\rm uloc}$ and  write that $$[ \dot M^{2,\frac 2 p q}, \dot W^{-\rho,\dot M^{r,\frac r p q}}]^{[\theta]}\subset [ \dot M^{2,\frac 2 p q}, \dot W^{-\rho,\dot M^{r,\frac r p q}}]_{\theta,\infty}\subset [ L^2_{\rm uloc}, \dot B^{-\rho-\frac {3p}{rq}}_{\infty,\infty}]_{\theta,\infty}$$

In order to finish the proof, we thus need to  use the machinery of energy control for suitable local Leray solutions \cite{LEM1, LEM2}. This will be done in the following sections, and we shall finish the proof in Section \ref{final} 
  $\diamond$
  
  \section{Weak local Leray solutions}\label{local}
  
  We recall basic results for local weak Leray solutions. We endow $L^2_{\rm uloc}$ with the norm
  $$ \|u\|_{L^2_{\rm uloc}}=\sup_{k\in \mathbb{Z}^3}  \|u \varphi_0(x-k)\|_2,$$ where $\varphi_0$ is a non-negative function in $\mathcal{D}$, suppported in a ball $B(0,R_0)$ and such that $\sum_{k\in\mathbb{Z}^3} \varphi_0(x-k)=1$.
 
  When $\vec u_0\in L^2_{\rm uloc}$,                                                                                                                                                                                                                                                                                                                                                                                                                                                                  
  proof of   existence of solutions   for the Navier--Stokes equations is based on mollification, energy estimates and compactness arguments (for details, see \cite{LEM2}, section 14.1) :
 \begin{itemize}
 \item we solve$$\partial_t \vec u_\epsilon+(\varphi_\epsilon*\vec u_\epsilon).\vec\nabla \vec u_\epsilon=\Delta\vec u_\epsilon -\vec\nabla p_\epsilon$$ with $\text{ div }\vec u_\epsilon=0$ and $\vec u_\epsilon(0,.)=\vec u_0$. Here, $\varphi\in \mathcal{D}$, $\int\varphi\, dx=1$ and $\varphi_\epsilon(x)=\frac 1{\epsilon^3} \varphi(\frac x \epsilon)$. Here $\vec\nabla p_\epsilon$ is given by the Leray projection : $$\vec \nabla p_\epsilon=- (\varphi_\epsilon*\vec u_\epsilon).\vec \nabla u_\epsilon +\mathbb{P} {\rm div }\left((\varphi_\epsilon*\vec u_\epsilon)\otimes\vec u_\epsilon\right).$$
  \item  the solution holds at least on an interval $(0,T_\epsilon)$ where $T_\epsilon$ depends on $\epsilon$ and on $\|\vec u_0\|_{L^2_{\rm uloc}}$ ($T_\epsilon=\min(1,C_0\frac {\epsilon^{3/2}}{\|\vec u_0\|_{L^2_{\rm uloc}}^2} $). Moreover,  we have the inequalities, for $k\in\mathbb{Z}^3$,
  \begin{equation*}\begin{split}
  \int \varphi_0(x-k) &\vert\vec u_\epsilon(t,x)\vert^2\, dx +2\int_0^t \int \varphi_0(x-k) \vert \vec\nabla\otimes\vec u_\epsilon(s,x)\vert^2\, dx\, ds
  \\ \leq &  \int \varphi_0(x-k) \vert\vec u_0(t,x)\vert^2\, dx  + C_1 \int_0^t\int_{\vert x-k\vert\leq R_0} \vert\vec u_\epsilon(s,x)\vert^2\, dx\, ds\\ &+ C_2 \int_0^t\int_{\vert x-k\vert>5R_0} \frac 1{\vert x-k\vert^4} \vert\vec u_\epsilon(s,x)\vert^3\, dx\, ds \\&+  C_2 \int_0^t\int_{\vert x-k\vert<5R_0}  \vert\vec u_\epsilon(s,x)\vert^3\, dx\, ds 
  \end{split}\end{equation*}
  \item defining
  $$ \alpha_\epsilon(t)=\|\vec u_\epsilon(t,.)\|_{L^2_{\rm uloc}}$$ and
  $$ \beta_\epsilon(t)=\sup_{k\in\mathbb{Z}^3} \left( \int_0^t \int \varphi_0(x-k) \vert\vec\nabla\otimes\vec u_\epsilon(s,x)\vert^2\, dx\, ds\right)^{1/2},$$
we get the inequality   \begin{equation*}\begin{split}
  \int \varphi_0(x-k) &\vert\vec u_\epsilon(t,x)\vert^2\, dx +\int_0^t \int \varphi_0(x-k) \vert \vec\nabla\otimes\vec u_\epsilon(s,x)\vert^2\, dx\, ds
  \\ \leq & \ \alpha_\epsilon(0)^2+ \frac 1 2 \beta_\epsilon(t)^2 + C_3 \int_0^t \alpha_\epsilon(s)^2\, ds +C_3 \int_0^t \alpha_\epsilon(s)^6\, ds
  \end{split}\end{equation*}
  so that
  $$ \beta_\epsilon(t)^2\leq 2 \alpha_\epsilon(0)^2+ 2C_3  \int_0^t \alpha_\epsilon(s)^2\, ds + 2C_3 \int_0^t \alpha_\epsilon(s)^6\, ds$$ and finally
    $$ \alpha_\epsilon(t)^2\leq 2 \alpha_\epsilon(0)^2+ 2C_3  \int_0^t \alpha_\epsilon(s)^2\, ds + 2C_3 \int_0^t \alpha_\epsilon(s)^6\, ds$$ Thus, as long as $8C_3 t<1$ and $128\, C_3 t  \|\vec u_0\|_{L^2_{\rm uloc}}^4<1$, we find that $\alpha_\epsilon(t)\leq 2 \|\vec u_0\|_{L^2_{\rm uloc}}$ and $\beta_\epsilon(t)\leq 2 \|\vec u_0\|_{L^2_{\rm uloc}}$.
\item the solution is then defined on $(0,\min(\frac 1{8C_3}, \frac 1{128\, C_3\, \|\vec u_0\|_{L^2_{\rm loc}}^4})$ and controlled independently from $\epsilon$. By Rellich theorem, wre find a subsequence that converges strongly in $(L^2_t L^2_x)_{loc}$ to a  suitable local Leray solution $\vec u$
 \end{itemize}

An important point is the following one : assume moreover that $\vec u_0\in {\rm BMO}^{-1}_0$ and that $\|e^{t\Delta}\vec u_0\|_{X_T}<\frac 1{4C_0}$ (where $C_0$ is the constant  of Theorem \ref{theoKT})
 then $\vec u_\epsilon$ is defined at least on $(0,T)$ and $\|\vec u_\epsilon\|_{X_T}\leq 2 \|e^{t\Delta}\vec u_0\|_{X_T}$. As $T$ does not depend on $\epsilon$, we see that the local Leray solution $\vec u$ satisfies $\vec u\in X_S$ with $S=\min(T, \frac 1{8C_3}, \frac 1{128\, C_3\, \|\vec u_0\|_{L^2_{\rm loc}}^4})$ and   $\|\vec u\|_{X_S}\leq 2 \|e^{t\Delta}\vec u_0\|_{X_T}$.

Similarly, if we assume   that $\vec u_0\in  \dot B^{-\gamma}_{\infty,\infty}$ with $-1<-\gamma<0$,    then $\vec u_\epsilon$ is defined at least on $(0,T)$ where $T= C \|\vec u_0\|_{\dot B^{-\gamma}_{\infty,\infty}}^{\frac 2 {1-\gamma}} $ and  $\sup_{0<t<T} t^{\frac \gamma 2} \|\vec u(t,.)\|_\infty \leq 2 \sup_{0<t<T} t^{\frac\gamma 2}\|e^{t\Delta}\vec u_0\|_{\infty}$. As $T$ does not depend on $\epsilon$, we see that the local Leray solution $\vec u$ satisfies the inequality  $\sup_{0<t<S} t^{\frac \gamma 2} \|\vec u(t,.)\|_\infty<+\infty$ where  $S=\min(T, \frac 1{8C_3}, \frac 1{128\, C_3\, \|\vec u_0\|_{L^2_{\rm loc}}^4})$.

\section{Comparison of local weak Leray solutions}
If $\vec u$ and $\vec v$ are two local weak Leray solutions, on $(0,T)$ with initial values $\vec u_0$ and $\vec v_0$, we would like to be able to estimate $\vec u(t,.)-\vec  v(t,.)$  from the estimation of $\vec u_0-\vec v_0$. 
This can be done only when at least one of the solutions is regular enough. We shall assume that $\vec u\in L^2_t L^\infty_x$. We sketch the computations described in \cite{LEM2}, section 14.4.

Define $\vec w=\vec u-\vec v$, $$ \alpha(t)=\|\vec w(t,.)\|_{L^2_{\rm uloc}}$$ and
  $$ \beta(t)=\sup_{k\in\mathbb{Z}^3} \left( \int_0^t \int \varphi_0(x-k) \vert\vec\nabla\otimes\vec w(s,x)\vert^2\, dx\, ds\right)^{1/2},$$ Using the suitability of $\vec v$ and the regularity of $\vec v$, we find
  (for $0<t<\min(1,T)$)
\begin{equation*}\begin{split}
  \int &\varphi_0(x-k) \vert\vec w(t,x)\vert^2\, dx +\int_0^t \int \varphi_0(x-k) \vert \vec\nabla\otimes\vec w (s,x)\vert^2\, dx\, ds
  \\ \leq & \ \alpha(0)^2+ \frac 1 2 \beta(t)^2 + C_1 \int_0^t \alpha(s)^2\, ds +C_2 \int_0^t \alpha(s)^6\, ds + C_3 \int_0^t \|\vec u(s,.)\|_\infty^2 \alpha(s)^2\, ds
  \end{split}\end{equation*}
where the constants $C_i$ do not depend on $T$, $\vec u$, nor on $\vec v$. Finally, we find
\begin{equation}\label{ineqweakstr}\begin{split} \alpha(t)^2\leq& 2 \alpha(0)^2 + 2 (C_1 +C_3 (\|\vec u\|_{L^\infty_t L^2_{\rm uloc}}+\|\vec v\|_{L^\infty_t L^2_{\rm uloc}})^4 )\int_0^t \alpha(s)^2\, ds\\&+ 2C_3 \int_0^t \|\vec u(s,.)\|_\infty^2 \alpha(s)^2\, ds.  \end{split}\end{equation}

We have the same estimate even if $\vec u$ is not intregrable near $t=0$.  Let us only assume that $\vec u\in L^2_tL^\infty_x$ on every $(\epsilon,T)$ with $\epsilon>0$. 
 Considering a time $t_0>0$ which a Lebesgue point for the functions $t\mapsto  \int \varphi(x-k)\vert\vec u(t,x)\vert^2\, dx$ and $t\mapsto  \int \varphi(x-k)\vert\vec u(t,x)\vert^2\, dx$, we have that $\vec u$ and $\vec v$ are local weak Leray solutions on $(t_0,T)$ with initial values $\vec u(t_0,.)$ and $\vec v(t_0,.)$. Thus, we shall find that, for $t>t_0$
 
\begin{equation}\label{ineqweakstr}\begin{split}   \int \varphi_0(x-k) \vert\vec w(t,x)\vert^2\, dx \leq& 2 \int\varphi_0(x-k) \vert\vec w(t_0,x)\vert^2\, dx\\&+ 2 (C_1 +C_3 (\|\vec u\|_{L^\infty_t L^2_{\rm uloc}}+\|\vec v\|_{L^\infty_t L^2_{\rm uloc}})^4 )\int_{t_0}^t \alpha(s)^2\, ds\\&+ 2C_3 \int_{t_0}^t \|\vec u(s,.)\|_\infty^2 \alpha(s)^2\, ds.  \end{split}\end{equation}

It is then enough to let $t_0$ go to $0$ and then take the supremum with respect to $k$.

  \section{Proof of Theorem \ref{main}} \label{final}

 We may now finish the proof. We consider two solutions of the Navier--Stokes equations with initial value $ \vec u_0 \in L^2\cap {\rm BMO}^{-1}_0\cap \dot B^{-s}_{\dot M^{1,q},\infty}$ with  $-s>-1+\frac 2 q$ (and $3<q<+\infty$) : we assume that  $\vec v$ is  a suitable   Leray solution and $\vec u$ is the mild solution in $X_T$.

 As $\vec v$ is suitable, $\vec v$ is a local Leray solution as well and we may estimate the $L^2_{\rm uloc}$ of $\vec u-\vec w$ :   defining $B(t)=\|\vec u(t,.)-\vec v(t,.)\|_{L^2_{\rm uloc}}^2$ and $$A(t)= 
 2 (C_1 +C_3 (\|\vec u\|_{L^\infty_t L^2_{\rm uloc}}+\|\vec v\|_{L^\infty_t L^2_{\rm uloc}})^4 ) + 2C_3  \|\vec u(t,.)\|_\infty^2,
$$
we get
$$ B(t)\leq \int_0^t A(s) B(s)\, ds. $$ As $\lim_{t\rightarrow 0} t A(t)=0$, we shall try to apply Lemma \ref{lemgron}. Thus, we shall use interpolation estimates to search for   a control of $B(t)$ as $O(t^{-\delta})$, in the spirit of Lemma  \ref{lembar}:.

 Recall that we have introduced the following numbers :
 \begin{itemize}
 \item $p$ such that $2<p<\frac{2q}3$
 \item $1-\theta$ the barycentric coordinate of $\frac 1 p$ in $[\frac 1 q, 1]$ :  $ \frac 1 p=(1-\theta) +\theta\frac 1 q$
  \item $-\sigma$ the corresponding point in $[-1 +\frac 2 q,-s]$ : $-\sigma=(1-\theta)(-s)+\theta(-1+\frac 2 q)$
 \item $-\gamma$ such that $-1+\frac 2 q<-\gamma<-\sigma$
 \item $r$ such that $p<r<+\infty$
 \item $1-\eta$ the barycentric coordinate of  $\frac 1 p$ in the segment  $[\frac 1 r,\frac 1 2]$ :
$ \frac 1 p= \frac {1-\eta}2+\frac \eta r$
\item $-\rho$ the corresponding point in $[-\gamma, 0]$ : $-\rho=\eta(-\gamma)$
\end{itemize}
 
 We have the following embeddings :
 \begin{itemize}
\item  $B^0_{2,2} \cap  {\rm BMO}^{-1}\subset   \dot B^{-1+\frac 2 q}_{q,q}$
\item $ \dot B^{-1+\frac 2 q}_{q,q}\cap \dot B^{-s}_{\dot M^{1,q},\infty} \subset  \dot B^{-\sigma}_{\dot M^{p,q},\infty}$  
\item $ \dot B^{-1+\frac 2 q}_{q,q}\cap  \dot B^{-\sigma}_{\dot M^{p,q},\infty}\ \subset \dot B^{-\gamma}_{\dot M^{p,q},1}\subset \dot W^{-\gamma,\dot M^{p,q}}.$
\item
$ W^{-\gamma,\dot M^{p,q}} =[ \dot M^{2,\frac 2 p q}, \dot W^{-\rho,\dot M^{r,\frac r p q}}]^{[\eta]}$ 
\item  $[ \dot M^{2,\frac 2 p q}, \dot W^{-\rho,\dot M^{r,\frac r p q}}]^{[\eta]}\subset [ \dot M^{2,\frac 2 p q}, \dot W^{-\rho,\dot M^{r,\frac r p q}}]_{\eta,\infty}\subset [ L^2_{\rm uloc}, \dot B^{-\rho-\frac {3p}{rq}}_{\infty,\infty}]_{\eta,\infty}$
\end{itemize}
If we take $r$ very large, we have $\eta=1-\frac 2 p+o(1)$ and $$-\rho-\frac{3p}{rq}=(1-\frac 2 p)(-\gamma) +o(1) \in ]-1,0[.$$

Thus far, we have seen that $\vec u_0\in   [ L^2_{\rm uloc}, \dot B^{-\lambda}_{\infty,\infty}]_{\eta,\infty}$ for some $\eta\in (0,1)$ and some $\lambda \in (0,1)$. We shall now estimate $\vec v-e^{t\Delta}\vec u_0$ when $\vec v$ is a weak local Leray solution.   For $0<\epsilon<1$, split $\vec u_0$ in $\vec \alpha_\epsilon+\vec\beta_\epsilon$ with
 $$ \|\vec\alpha_\epsilon\|_{L^2_{\rm uloc}}\leq C_0 \epsilon^{\eta}\text{ and } \|\vec \beta_\epsilon\|_{\dot B^{-\lambda}_{\infty,\infty}}\leq C_0 \epsilon^{\eta-1}$$ where $C_0$ does not depend on $\epsilon$ (but depends on $\vec u_0$).
 
 As $\beta_\epsilon=\vec u_0-\alpha_\epsilon$, we have $\| \vec \beta_\epsilon\|_{L^2_{\rm uloc}}\leq \| \vec u_0\|_{L^2_{\rm uloc}} + C_0$ (an estimation which does not depend on $\epsilon$) and we know that we have a (restricted) weak Leray solution of the Navier--Stokes equations $\vec U_\epsilon$ with initial value $\vec\beta_\epsilon$  such that $\sup_{0<t<T_0} \|\vec U_\epsilon(t,.)\|_{L^2_{\rm uloc}}\leq C_1$,  where $T_0$ and $C_1$ depends only on $\|\vec u_0\|_{L^2_{\rm uloc}}$ (and not on $\epsilon$).
 
 As $\beta_\epsilon\in \dot B^{-\lambda}_{\infty,\infty}$,, we have as well 
   that $\| \vec U_\epsilon(t,.)\|_\infty\leq C_2 t^{-\frac {\lambda }2}  \epsilon^{\eta-1}$ on an interval $(0,T_\epsilon)$  with $T_\epsilon^{\frac {1-\lambda}2} \epsilon^{\eta-1}= C_3$ It is then easy to check that $\|\vec U_\epsilon-e^{t\Delta}\vec\beta_\epsilon\|_{L^2_{\rm loc}} \leq C_4 t^{(1-\lambda)/2} \epsilon^{\eta-1}.$

Using our results on comparison of suitable local Leray solutions, we find that 
  we get the following inequality for  $\vec W_\epsilon=\vec v-\vec U_\epsilon$ and $A_\epsilon(t)=\sup_{k\in \mathbb{Z}^3} \int \varphi_0(x-k) \vert \vec W_\epsilon(t,.x)\vert^2\, dx$    : 
\begin{equation}\label{ineqweakstr}  A_\epsilon(t)^2\leq 2 \|\alpha_\epsilon\|_{L^2_{\rm uloc}}^2+ C_5\int_0^t A_\epsilon(s)^2\, ds +  C_6 \int_0^t \|\vec u(s,.)\|_\infty^2 A_\epsilon(s)^2\, ds.   \end{equation}

Thus, we get
$$    A_\epsilon(t) \leq C_7 \epsilon^{2\eta} + C_8\int_0^t A_\epsilon(s)^2\, ds +  C_9 \epsilon^{2(\eta-1)} \int_0^t s^{-\lambda} A_\epsilon(s)^2\, ds.   $$
  so that
$$    A_\epsilon(t) \leq C_7 \epsilon^{2\eta}  e^{C_8 t}  e^{ C_9 \epsilon^{2(\eta-1)}\frac{t^{1-\lambda}}{1-\lambda}}.
$$

   Now, for $\tau<1$, take $\epsilon=\tau^\mu$ with $\frac{1-\lambda} 2+\mu(\eta-1)>0$. We find that, for $0<t< \min(T_\epsilon,T_0)$ with 
 $$ T_\epsilon=  C_3^{\frac 2{1-\lambda}}\epsilon^{\frac{2 (1-\eta)}{1-\lambda}}= C_3^{\frac 2{1-\lambda}} \tau^{\mu \frac{2 (1-\eta)}{1-\lambda}}\  [\text{ where }  \mu \frac{2 (1-\eta)}{1-\lambda}<1],$$ we have the inequality
\begin{equation*}\begin{split} \|\vec v-e^{t\Delta}\vec u_0\|_{L^2_{\rm uloc}}\leq & \|\alpha_\epsilon\|_2+\|\vec W_\epsilon\|_2+ \|\vec U_\epsilon-e^{t\Delta} \vec \beta_\epsilon\|_2
\\ \leq &  C_{10} \epsilon^\eta (1+e^{  C_{11} \frac {t^{1-\lambda}}{T_\epsilon^{1-\lambda}}}) +C_4 t^{\frac{1-\lambda}2} \epsilon^{\eta-1}
\end{split}\end{equation*}
 If $\tau$ is small enough, we have $\tau<T_\epsilon$ and we find
 $$ \|\vec v(\tau,.)-e^{\tau\Delta}\vec u_0\|_2\leq    C_0\tau^{\mu\eta} (1+e^{ \frac {C_1^2}{2(1-\lambda }} )+C_4 \tau^{\frac{1-\lambda}2+\mu(\eta-1)}$$
 The theorem  is proved.


\end{document}